# Subjective-Objective Median-based Importance Technique (SOMIT) to Aid Multi-Criteria Renewable Energy Evaluation

Ding Ding [a, †, 0000-0002-1038-7558], Yang Li [a, †, 0009-0007-4368-9732], Poh Ling Neo [a, †, 0000-0003-3710-5976], Zhiyuan Wang [a, †, *, 0000-0001-7867-7626], Chongwu Xia [a, †, 0000-0002-9343-4670]

[a] School of Business, Singapore University of Social Sciences, Singapore 599494, Singapore

\* Corresponding author: Zhiyuan Wang, Email: zywang@suss.edu.sg

[†] All authors contributed equally to this work and are listed in alphabetical order by last name

## Abstract

Accelerating the renewable energy transition requires informed decision-making that accounts for the diverse financial, technical, environmental, and social trade-offs across different renewable energy technologies. A critical step in this multi-criteria decision-making (MCDM) process is the determination of appropriate criteria weights. However, deriving these weights often solely involves either subjective assessment from decision-makers or objective weighting methods, each of which has limitations in terms of cognitive burden, potential bias, and insufficient contextual relevance. This study proposes the subjective-objective median-based importance technique (SOMIT), a novel hybrid approach for determining criteria weights in MCDM. By tailoring SOMIT to renewable energy evaluation, the method directly supports applied energy system planning, policy analysis, and technology prioritization under carbon neutrality goals. The practical utility of SOMIT is demonstrated through two MCDM case studies on renewable energy decision-making in India and Saudi Arabia. Using the derived weights from SOMIT, the technique for order preference by similarity to ideal solution (TOPSIS) method ranks the renewable energy alternatives, with solar power achieving the highest performance scores in both cases (e.g., 0.5725 in the India case). The main contributions of this work are five-fold: 1) the proposed SOMIT reduces the number of required subjective comparisons from the conventional quadratic order to a linear order; 2) SOMIT is more robust to outliers in the alternatives-criteria matrix (ACM); 3) SOMIT balances subjective expert knowledge with objective data-driven insights, thereby mitigating bias; 4) SOMIT is inherently modular, allowing both its individual parts and the complete approach to be seamlessly coupled with a wide range of MCDM methods commonly applied in energy systems and policy analysis; 5) a dedicated Python library, pysomit, is developed for SOMIT, providing an accessible and efficient tool to implement SOMIT in practical renewable energy evaluation and decision-support applications.

**Keywords:**

Renewable Energy Evaluation

Multi-Criteria Decision-Making

Subjective-Objective Median-based Importance Technique

TOPSIS Method

Carbon Neutrality

**1. Introduction**

In the face of escalating climate change, dwindling fossil fuel reserves, and rising energy demands, renewable energy has emerged as a cornerstone of global efforts to achieve a sustainable and resilient energy future (Halkos & Gkampoura, 2020). Unlike conventional energy sources, renewables (e.g., solar, wind, hydropower, geothermal, and biomass) offer the dual advantage of reducing greenhouse gas (GHG) emissions and enhancing energy security (Tutak & Brodny, 2022). Over the past decade, significant technological advancements, policy incentives, and cost reductions have accelerated the deployment of renewable energy systems worldwide (Hassan et al. 2024). As countries strive to meet their carbon neutrality goals and diversify their energy portfolios, renewable energy is no longer viewed merely as an alternative, but as a central pillar of modern energy strategies (Rabbi et al., 2022).

International collaborations, such as the Paris Agreement, along with national policies and investment flows, have further solidified the role of renewables in shaping long-term energy planning. Moreover, increasing public awareness and demand for cleaner energy have pushed both governments and industries to fast-track the transition toward low-carbon energy solutions. This global trend is particularly prominent in Asia, home to over half of the world's population and some of the fastest-growing economies. Countries like China, India, and members of the ASEAN (Association of Southeast Asian Nations) bloc are experiencing surging energy demand due to rapid urbanization, industrialization, and rising living standards (Anwar et al., 2022). At the same time, the region holds vast renewable energy potential, which makes the transition to clean energy both a necessity and an opportunity for sustainable development (Li et al., 2020). These challenges underscore the need for advanced decision-support tools tailored to energy planning and policy. In particular, determining appropriate weights across environmental, technical, financial, and social criteria is a critical step in aligning renewable energy strategies with carbon neutrality and sustainability objectives.

Despite the growing momentum behind renewable energy, the selection and prioritization of specific technologies remain a complex and multifaceted challenge (Sitorus & Brito-Parada, 2020). Decision-makers must account for a wide range of criteria that often conflict with one another, for instance, minimizing financial costs while maximizing environmental sustainability (Terlouw et al., 2023) or balancing technical performance with social acceptance (Hassan et al., 2022). These criteria usually span several categories and domains, including financial factors (e.g., capital cost, operating expenses), technical attributes (e.g., efficiency, technology maturity), environmental impacts (e.g., GHG emissions, total land use), and social considerations (e.g., job creation, social acceptance). The inherent trade-offs among these dimensions make the evaluation process far from straightforward, especially in contexts where stakeholder priorities and regional conditions vary significantly. As such, selecting the most suitable renewable energy option requires a structured and transparent approach that can accommodate both objective quantitative data and subjective qualitative judgments (Karunathilake et al., 2019).

In this regard, multi-criteria decision-making (MCDM) methods have been widely adopted worldwide for renewable energy evaluation due to their ability to handle complex and multiple criteria for informed decision analysis. Through an extensive review of the literature on MCDM methods applied to renewable energy source evaluation, we found that determining the weights of evaluation criteria remains one of the most critical and challenging tasks for decision-makers. These weights directly influence the outcome of the decision process, yet assigning them often involves subjectivity, inconsistency, bias, or a heavy cognitive burden. Moreover, it is observed that existing MCDM applications often rely solely on either subjective or objective weighting methods, both of which may introduce bias or overlook critical aspects of the decision context.

To address these shortcomings, this study proposes a novel hybrid weighting method, referred to as the Subjective-Objective Median-based Importance Technique (SOMIT), which integrates subjective weights derived from decision-maker's assessments with objective data-driven weights obtained from the ACM through a three-part framework. In Part I, decision-makers mainly perform pairwise comparisons of all criteria against a single median-importance reference criterion, significantly reducing the number of comparisons to a linear order (i.e., on the order of magnitude of $n$ comparisons, that is, a linear scale). In Part II, objective weights are derived from the normalized median-based dispersion of each criterion in the ACM, capturing data-driven variation and importance. In Part III, subjective and objective weights are combined using a multiplicative synthesis, producing the final balanced weights for all criteria. The main advantages of SOMIT are: 1) it minimizes cognitive load by reducing subjective pairwise comparisons from quadratic to linear complexity; 2) it is more robust to outliers in the ACM; 3) it balances human expertise with empirical data, mitigating the possible bias of purely subjective or objective weighting methods; 4) it features a modular design, enabling either each part or the complete approach to be seamlessly coupled with a wide range of MCDM methods; 5) an easy-to-use Python library, pysomit, has been developed to facilitate the practical use of the SOMIT method.

As aforementioned, while earlier studies in the literature have demonstrated the value of MCDM in renewable energy evaluation, few have explicitly addressed the weighting stage as an energy policy challenge. In this context, the proposed SOMIT method not only advances methodological rigor but also enhances applied energy decision-making by reducing the cognitive burden of experts and improving robustness, while ensuring that critical energy system dimensions (e.g., GHG mitigation, land use, social acceptance) are incorporated. Notably, this work positions SOMIT in the applied energy context, where financial, technical, environmental, and social trade-offs need to be balanced for achieving carbon neutrality and ensuring sustainable energy transitions. Through the renewable energy case studies in India and Saudi Arabia, this work showcases SOMIT's direct applicability in renewable energy planning and decision-making. In this way, it explicitly connects the methodological contribution to the applied aspects of energy by supporting optimal use of energy resources through a structured decision-making framework that accounts for multidimensional impacts of renewables, aligning with global shifts toward a net zero future.

The rest of this article is structured as follows. Section 2 reviews the relevant literature on MCDM and weighting methods in renewable energy applications. Section 3 introduces the proposed SOMIT weighting method, describing its three-part approach for deriving criteria weights. This section also briefly outlines the widely used technique for order preference by similarity to ideal solution (TOPSIS) MCDM method, which is employed in this work to rank renewable energy alternatives using the weights generated by SOMIT. Section 4 applies the methods to two case studies in India and Saudi Arabia, presenting the detailed results and drawing out key implications for renewable energy planning. Section 5 discusses the findings in greater depth and compares SOMIT with two established weighting methods. Finally, Section 6 concludes the current study and provides recommendations for future research.

## 2. Literature Review

Renewable energy plays a central role by replacing fossil fuels with cleaner options (Tian et al., 2022), in order to achieve the goal of carbon neutrality that requires reducing emissions at their source. Carbon neutrality refers to balancing the amount of GHG released into the atmosphere with an equivalent amount removed or offset, resulting in a net zero future (Chen et al., 2024). These renewables generate energy with little to no carbon output and are naturally replenished. Many countries around the world are committing to carbon neutrality as part of their strategies for a sustainable future. Nations such as the European Union (EU) members, Japan, South Korea, Canada, and China have set the mid- to long-term targets, often aiming for net zero emissions by 2050 or 2060 (Dafnomilis et al., 2023). They are also advancing renewable energy and green technologies to align with these goals. These commitments involve transitioning away from fossil fuels, investing in clean energy, improving energy efficiency, and protecting natural carbon sinks like forests. By pursuing carbon neutrality, countries not only aim to mitigate the effects of climate change but also to promote innovation, enhance energy security, and foster green economic growth, shaping a more resilient and sustainable future for the planet (Qamruzzaman & Karim, 2024).

For MCDM applications in renewable energy systems, Zhang et al. (2019) applied a hybrid MCDM approach to evaluate renewable energy projects in Fujian, China, where the assessment focused on criteria spanning economic (e.g., investment cost), technical (e.g., technical maturity and efficiency), environmental (e.g., GHG emissions), and social (e.g., job creation). Rahim et al. (2020) used MCDM methods to select the most suitable renewable energy source in Malaysia based on criteria including technical, economic, environmental, and social aspects. Solangi et al. (2021) used an integrated MCDM techniques to assess barriers and prioritize strategies for renewable energy deployment in Pakistan. Bilgili et al. (2022) applied MCDM methods to evaluate and rank renewable energy options for sustainable development in Turkey, considering many criteria across economic, technical, environmental, social, and political dimensions. Prayogo el al. (2023) applied MCDM methods in the enhancement of solar photovoltaic and wind turbine systems in Indonesia. Ridha et al. (2023) utilized MCDM methods to evaluate and rank Pareto-optimal front solutions for hybrid photovoltaic (PV), wind, and battery systems in Malaysia and South Africa, considering multiple technical, economic, environmental, and reliability criteria. Amiri et al. (2024) employed MCDM methods for renewable energy source selection in Saudi Arabia, taking into account criteria such as cost-effectiveness, land conservation, and design simplicity. Li et al.

(2024) applied the TOPSIS method for the selection of wind turbines based on efficiency and economic performance for wind energy development. Vaziri Rad et al. (2024) adopted a hybrid MCDM framework to evaluate and prioritize sustainable energy system configurations for agri-photovoltaic applications, considering environmental, economic, technical, energy security, and social objectives. A recent review article by Kumar and Pamucar (2025) highlighted the growing attention on sustainable energy planning using MCDM approaches. Sahoo et al. (2025) also conducted a comprehensive review demonstrating the global application of MCDM methods in planning renewable energy and addressing sustainability challenges. Tabrizi et al. (2025) utilized MCDM methods to evaluate renewable energy adoption across G7 countries (i.e., the seven major developed economies: Canada, France, Germany, Italy, Japan, the United Kingdom, and the United States). Murtaja et al. (2025) applied MCDM methods to select optimal renewable energy sources in India based on six criteria including rural life upliftment, job creation, use of local resources, social acceptance, and pollution/noise impacts.

When deriving the weights of criteria, the Analytic Hierarchy Process (AHP) method, proposed by Saaty (1990), is the most popular subjective weighting method according to a review by Hafezalkotob et al. (2019). AHP requires decision-makers to conduct pairwise comparisons to express the relative importance of criteria. However, AHP presents several notable limitations. First, it requires a full set of pairwise comparisons across all criteria (on the order of magnitude of $n^2$ comparisons for the $n$ criteria, that is, a quadratic scale), which becomes increasingly cumbersome when the number of criteria grows (Carmone et al., 1997). This imposes a significant cognitive load on decision-makers. Second, using the original AHP (Saaty, 1990), deriving weights from the pairwise comparison matrix involves solving an eigenvalue–eigenvector problem, which is computationally nontrivial and practically inaccessible for most decision-makers without the aid of computational tools, especially when the matrix exceeds a manageable size such as 3×3. While approximation techniques for AHP exist (Wang et al., 2020), they may yield weights that deviate from those obtained via the original eigenvalue method, potentially affecting decision accuracy. In contrast, objective weighting methods, such as the widely used criteria importance through intercriteria correlation (CRITIC) (Brodny & Tutak, 2023; Wang et al., 2022), bypass human judgment and derive weights directly from the alternatives-criteria matrix (ACM). As purely data-driven approaches, they do not incorporate expert knowledge or contextual priorities, potentially overlooking valuable subjective insights.

## 3. Methodology

The methodology in this study follows a structured and adaptable decision analysis framework, as illustrated in the flowchart in Figure 1. This framework is applicable to many decision-making problems that require the evaluation and ranking of alternatives across multiple criteria.

In this work, the framework is specifically adapted for the evaluation of renewable energy, as shown in Figure 1. The process begins with the construction of the ACM, which captures the performance of each alternative (e.g., solar and wind) with respect to the relevant energy-specific criteria (e.g., technical and environmental) and its sub-criteria. This ACM forms the foundation for analysis. Then, the proposed SOMIT method is employed to determine the relative importance (i.e., weight) of each criterion. SOMIT, by default, is a hybrid approach that combines both sets of subjective and objective weights, although either set can also be applied independently within the framework. SOMIT reduces cognitive burden by requiring only a linear number of pairwise comparisons and enhances reliability by incorporating data-driven insight. The result is a set of balanced and context-sensitive weights that reflect both expert judgment and empirical evidence. Subsequently, these criteria weights are integrated into one or more MCDM methods, such as the TOPSIS introduced in Hwang and Yoon (1981), the MARCOS (measurement of alternatives and ranking according to compromise solution) proposed by Stević et al. (2020), and the PROBID (preference ranking on the basis of ideal-average distance) created by Wang et al. (2021), to evaluate the alternatives and assist in providing insights for renewable energy planning and policy. This produces a transparent recommendation for the top-performing alternative(s).

The framework is modular, interpretable, and efficient. By combining expert knowledge with structured data analysis, the methodology ensures a comprehensive and robust basis for complex MCDM problems. Unlike conventional methods that rely solely on expert opinion (e.g., AHP) or statistical dispersion within ACM itself (e.g., CRITIC), SOMIT offers a balanced integration. Additionally, the use of the median enhances robustness against outliers, ensuring that final weights more accurately represent consensus views while remaining grounded in the objective ACM data. This is particularly valuable in renewable energy planning, where evaluation criteria often span conflicting perspectives such as cost, GHG emissions, land usage, and social acceptance. Although SOMIT is generalizable to multiple domains, in this study it is explicitly tailored to renewable energy evaluation, ensuring that the weighting process reflects the applied energy concerns of carbon neutrality, energy planning, and sustainable technology adoption.

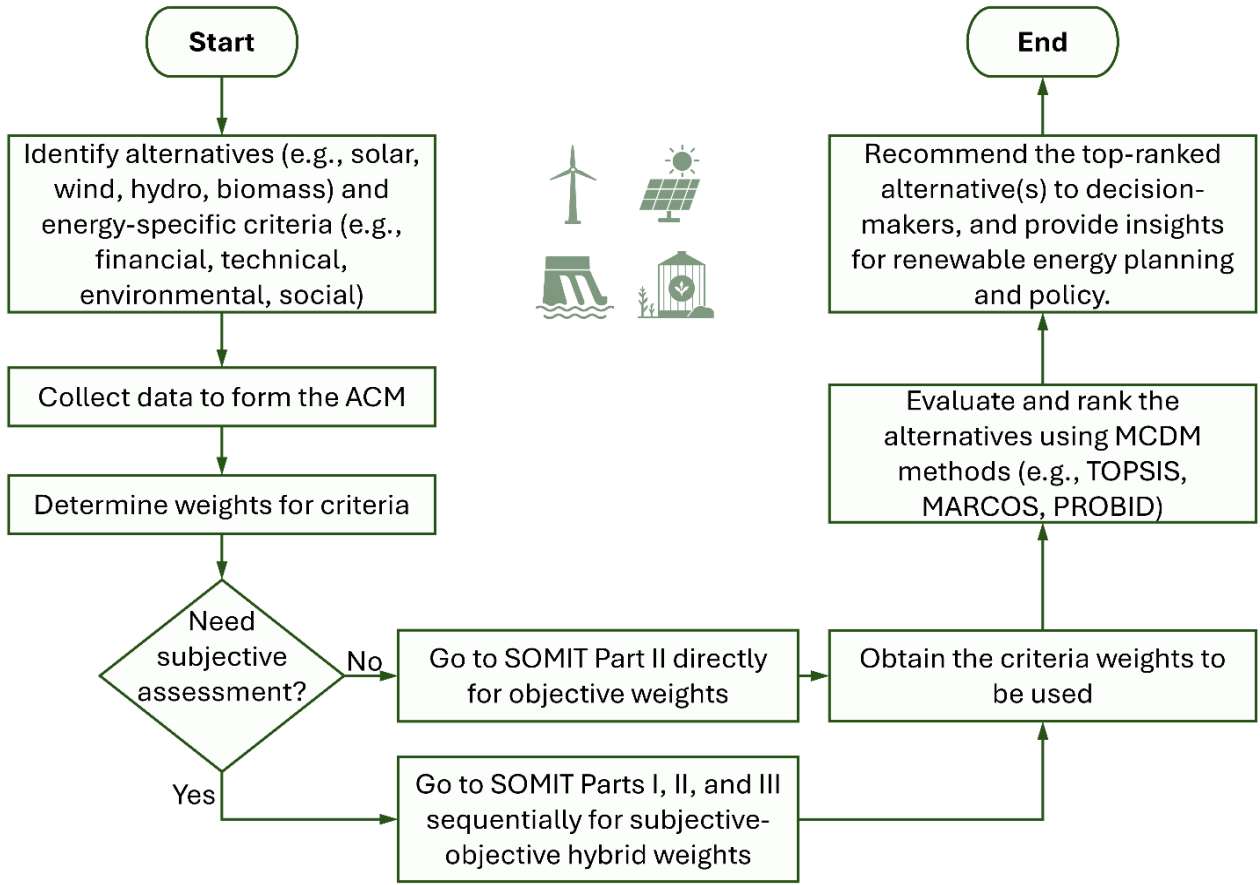

Figure 1. Flowchart of the decision analysis framework integrating ACM construction, SOMIT weighting, and MCDM ranking for renewable energy evaluation.

### 3.1. Subjective-Objective Median-based Importance Technique (SOMIT)

This study proposes a novel hybrid method for determining the weights of criteria, SOMIT, which is structured into three parts. In the first part, the decision-maker provides subjective assessments of the relative importance of criteria using a 1–9 comparison scale to derive the subjective weights, reflecting expert judgment. This part requires only $n$ comparisons, where $n$ is the number of criteria. The second part involves calculating objective weights directly from the alternatives-criteria matrix (ACM). This is done without any direct involvement from the decision-maker, relying solely on the performance data of alternatives across the criteria. In the third part, the subjective and objective weights are integrated using a multiplicative synthesis approach. This allows for a balanced representation of expert judgment and data-driven evidence, producing the final weights for the criteria. The main advantages of SOMIT are that it reduces the cognitive burden on the decision-maker by requiring only $n$ pairwise comparisons, as opposed to the conventional methods that typically involve an order of magnitude of $n^2$ comparisons to build a full pairwise matrices. SOMIT is more robust to outliers in ACM and promotes consistency and balance by integrating both subjective judgments from the decision-maker and objective data-driven insights from the ACM, thereby mitigating the risk of bias associated with using only either subjective or objective weights. SOMIT is easy to couple with a wide range of MCDM methods in various applications. An open-source Python library, pysomit, was developed for the practical use of this method. The three parts of SOMIT unfold as follows.

**SOMIT Part I: Derive subjective weights based solely on decision-maker's assessment**

Suppose a decision-maker is evaluating a total of $n$ criteria to determine the importance (i.e., weight) of each criterion subjectively. It is worth mentioning that the SOMIT Part I relies solely on the decision-maker's assessment of the criteria themselves, without considering the alternatives.

**Step 1:** The decision-maker first identifies one criterion that he/she subjectively perceives as being of median or moderate level of importance (i.e., neither highly important nor unimportant) relative to the others. This median criterion is denoted as $C_d$, where $d \in \{1,2,\ldots,n\}$.

**Step 2:** Every other criterion ($C_i$, where $i \in \{1,2,\ldots,n\}$ and $i \neq d$) is then compared with the median criterion ($C_d$) using the widely adopted pairwise comparison 9-point scale, as shown in Table 1. Each comparison reflects the relative importance of criterion $C_i$ with respect to $C_d$. This comparison value is represented as $a_{id}$.

For instance, if the decision-maker judges a criterion $C_i$ is moderately more important than $C_d$, then the comparison value is $a_{id} = 3$. This can also be interpreted as $C_i$ is three times as important as $C_d$; reciprocally, $C_d$ is deemed one-third as important as $C_i$, i.e., $a_{di} = 1/3$.

It is important to note that, since all comparisons in this step are made relative to the median criterion, it is generally expected that the decision-maker would rate approximately half of the criteria as more important than $C_d$ (i.e., $a_{id} \geq 1$, within the range [1, 9]), and the other half as less important (i.e., $a_{id} \leq 1$, within the range [1/9, 1]). However, this "half–half" distribution is viewed as a soft guideline rather than a strict constraint in this step. Some degree of deviation or inconsistency is acceptable, given the subjective nature of human judgment in this process.

Table 1. Pairwise comparison scale.

| Scale | Compare between criteria $C_i$ and $C_j$ |
|---|---|
| 1 | $C_i$ is **equally** important to $C_j$ |
| 3 | $C_i$ is moderately **more** important than $C_j$ |
| 5 | $C_i$ is strongly **more** important than $C_j$ |
| 7 | $C_i$ is very strongly **more** important than $C_j$ |
| 9 | $C_i$ is extremely **more** important than $C_j$ |
| 1/3 | $C_i$ is moderately **less** important than $C_j$ |
| 1/5 | $C_i$ is strongly **less** important than $C_j$ |
| 1/7 | $C_i$ is very strongly **less** important than $C_j$ |
| 1/9 | $C_i$ is extremely **less** important than $C_j$ |
| 2, 4, 6, 8, 1/2, 1/4, 1/6, 1/8 | Allowed; Denoting intermediate value between adjacent scales shown above |
| Any decimal value within the range [1/9, 9], e.g., 0.21 and 3.5 | Allowed; Denoting finer comparison scales |

**Step 3:** Identify the criteria with the highest and lowest comparison values relative to $C_d$ from Step 2, denoted as $C_h$ (highest) and $C_l$ (lowest), respectively, where $h, l \in \{1,2,...,n\}$. The decision-maker then directly assesses the relative importance of $C_h$ compared to $C_l$ using the same comparison scale provided in Table 1, denoted as $a_{hl}$.

We deliberately single out the $C_h$ and $C_l$ criteria, requiring the decision-maker to make one additional direct comparison between them. This is mainly because these two criteria sit at the extreme ends of the subjective weight distribution, and the subsequent MCDM process is particularly more sensitive to them. The relative importance of these two criteria effectively bounds the feasible range of all other weights and therefore functions as a "gate" for the entire subjective weighting scheme. By asking the decision-maker to directly compare them and state the $a_{hl}$ value, we create strong upper- and lower-bound anchors. If these anchors conflict sharply with the indirect comparisons from Step 2, it will expose inconsistency in the following Step 4 for the constrained optimization.

**Step 4:** Ideally, for all $n$ criteria, with $i, j \in \{1,2,...,n\}$, the subjective weights of any two criteria ($w_i^s, w_j^s$) should be determined such that the pairwise comparison value $a_{ij}$ (i.e., the relative importance of criterion $C_i$ compared to $C_j$) equals exactly the ratio of their weights, $w_i^s/w_j^s$ (this ratio essentially also reflects how important $C_i$ is relative to $C_j$). Mathematically, this ideal scenario is expressed as: $a_{ij} = w_i^s/w_j^s$, or equivalently: $a_{ij}w_j^s = w_i^s$, where the superscript ($s$) in $w_i^s$ and $w_j^s$ indicates that these weights are derived from subjective assessment.

However, pairwise comparisons made by human decision-makers are inherently subjective and often contain a degree of vagueness or inconsistency. For instance, if a decision-maker considers $C_1$ to be two times more important than $C_2$, and $C_2$ to be four times more important than $C_3$, then logically, $C_1$ should be eight times more important than $C_3$. This logical condition is known as the transitive property, that is, $a_{ij}a_{jk} = a_{ik}$, for all $i,j,k \in \{1,2,\ldots,n\}$. Despite this ideal case, it is recognized that human judgment may not strictly adhere to this property due to cognitive biases, subjectivity, and/or inadvertent inconsistencies.

Therefore, the following constrained optimization problem is formulated in Step 4, aiming to minimize the deviation to the ideal scenario.

$$\text{Minimize} \quad z = \sum_{(i,j)\in \mathcal{A}} \left(a_{ij}w_j^s - w_i^s\right)^2 \tag{1}$$

where $\mathcal{A}$ is the set of all pairs $(i,j)$ for which $a_{ij}$ has been assigned a pairwise comparison value in Steps 2 or 3.

$$\text{Subject to} \quad \sum_{j=1}^{n} w_j^s = 1, \quad \text{where } 0 \leq w_j^s \leq 1 \tag{2}$$

This problem can be solved manually using the method of Lagrange multipliers (Bertsekas, 1982), which transforms it into a system of linear equations. Solving this system produces the solution, that is, a set of subjective weights $\{w_1^s, w_2^s, \ldots, w_n^s\}$.

The method of Lagrange multipliers tackles constrained optimization by converting it into an unconstrained one. To find the extrema of a scalar objective function $z$ subject to an equality constraint of $g = 0$, we introduce a single multiplier $\alpha$ and form the Lagrangian $L = z - \alpha g$. At any optimum that respects the constraint, infinitesimal movements in the underlying variables must leave $z$ unchanged to first order and keep $g = 0$ satisfied; algebraically, this is enforced by requiring that all first-order derivatives of the Lagrangian vanish simultaneously. In symbols,

$$\nabla z - \alpha \nabla g = \mathbf{0}, \quad g = 0 \tag{3}$$

Solving these conditions yields candidate points where the contours of $z$ are tangent to the constraint manifold defined by $g$, ensuring that any feasible deviation cannot increase or decrease $z$ to first order. The multiplier $\alpha$ also provides sensitivity information: it quantifies how much the optimal value of function $z$ would change per unit relaxation of the constraint.

Alternatively, the constrained optimization problem can also be solved with the aid of computational tools such as Microsoft Excel Solver, Python, MATLAB, or other software capable of handling constrained optimization problems. One advantage of SOMIT Part I is that it requires only $n$ comparisons from the decision-maker to determine the subjective weights of $n$ criteria, significantly reducing the cognitive burden as opposed to the conventional AHP approach.

**SOMIT Part II: Derive objective weights from alternatives-criteria matrix (ACM)**

Before proceeding to Part II of SOMIT, it is important to clarify the use of the term "objective weights". In this context, "objective" is an adjective to indicate that the weights are not influenced by personal or subjective judgment. It should not be confused with "objective" as a noun, as in the context of multi-objective optimization (MOO) problems (Qi et al., 2017; Su et al., 2023; Terada et al., 2025), where it typically means a goal or target. Once the ACM is established, these weights are derived purely from the data contained in the ACM, independent of individual preferences or opinions.

**Step 5:** Normalize the ACM with $m$ rows (alternatives) and $n$ columns (criteria) by applying the Max-Min normalization. Note that, throughout the Parts II and III of SOMIT, the row and column indices (i.e., $i$ and $j$) always satisfy that $i \in \{1,2,\dots,m\}$ and $j \in \{1,2,\dots,n\}$.

$$F_{ij} = \frac{f_{ij} - \min\limits_{k \in \{1,2,\dots,m\}} f_{kj}}{\max\limits_{k \in \{1,2,\dots,m\}} f_{kj} - \min\limits_{k \in \{1,2,\dots,m\}} f_{kj}} \quad (4)$$

where $f_{ij}$ is the value of the $i^{th}$ alternative under the $j^{th}$ criterion (i.e., $C_j$) in the original ACM, and $F_{ij}$ is the corresponding normalized value in the normalized ACM. Max-Min normalization can preserve the relative dispersion of each criterion while transforming values of all criteria onto a common scale between 0 and 1.

**Step 6:** Find the median value of each normalized criterion (i.e., each column in the normalized ACM), and then compute the average absolute deviation from the median (AADM) for each criterion. This AADM value, denoted as $r_j$, quantifies the spread of criterion $C_j$.

$$r_j = \frac{1}{m}\sum_{i=1}^{m}|F_{ij} - \tilde{F}_j| \quad (5)$$

where $\tilde{F}_j = \text{median}(F_{1j}, F_{2j}, \ldots, F_{mj})$ is the median value of the normalized $C_j$ across all $m$ alternatives.

**Step 7:** Compute the objective weight $w_j^o$ for each criterion by normalizing its $r_j$ value, where the superscript ($o$) in $w_j^o$ indicates that these weights are derived objectively from ACM.

$$w_j^o = \frac{r_j}{\sum_{k=1}^{n} r_k} \tag{6}$$

In general, a criterion that exhibits greater spread across alternatives is considered to carry more information and is therefore assigned a higher weight. The underlying rationale is that greater dispersion in a criterion means it can better distinguish among the alternatives. In other words, criterion values that are closely clustered offer limited insight for differentiation, whereas values that span a broader range reveal clearer performance distinctions. After this step, the set of objective weights is obtained: $\{w_1^o, w_2^o, \ldots, w_n^o\}$.

Here, in Part II of SOMIT, we use the median as the reference level when computing the AADM value of each criterion because median is a robust measure of central tendency, which is unaffected by extreme outliers or skewed distributions. By anchoring spread on this resistant center, it reflects the typical dispersion of a criterion without being inflated by a few aberrant alternatives, producing weights that are more stable and representative of the bulk of the data. This robustness is especially useful in multi-criteria settings where some alternatives may have atypically high or low values, and it aligns with the decision-analytic goal of preventing a single anomalous entry from distorting the overall weighting scheme.

**SOMIT Part III: Integrate subjective and objective weights to derive the final weights**

**Step 8:** In the final step, the subjective and objective weights are integrated to compute the final weight $w_j$ for each criterion as follows:

$$w_j = \frac{w_j^s \cdot w_j^o}{\sum_{k=1}^{n} (w_k^s \cdot w_k^o)} \tag{7}$$

It is worth noting that SOMIT is purposefully designed to be modular, allowing decision-makers the flexibility to use either the subjective weights from Part I, the objective weights from Part II, or the combined weights from Part III in their subsequent MCDM process. While this choice is discretionary, using the combined weights is generally recommended, as it incorporates both the decision-makers' subjective assessments and the objective data-driven insights provided by the ACM.

To facilitate the application of the SOMIT method, we have developed a user-friendly Python library, named pysomit, which is available through the Python Package Index

([https://pypi.org/](https://pypi.org/)). Users can easily install the library using the standard pip command (pip install pysomit), as shown in Figure 2.

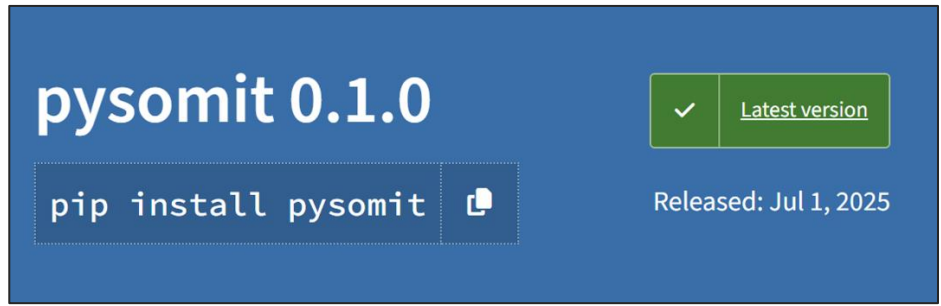

Figure 2. Screenshot of the pysomit library (version 0.1.0) hosted on the standard Python Package Index (PyPI) repository ([https://pypi.org/project/pysomit/](https://pypi.org/project/pysomit/)) and the installation command using pip.

For decision-makers who are unfamiliar with Python, one of the most accessible ways to start using this powerful tool is via Google Colaboratory ([https://colab.google/](https://colab.google/)). By opening a new notebook there, users can easily execute Python code to install and run functions from the pysomit library. As shown in Figure 3, pysomit provides an interactive prompt that guides users through the process of supplying subjective assessments required for SOMIT Part I.

```
1 # Uncomment next line if you haven't installed the pysomit library
2 # !pip install pysomit
3
4 from pysomit import somit_i, somit_ii, somit_iii
5
6 # SOMIT Part I
7 subjective_weights = somit_i()
8
9 # SOMIT Part II
10 objective_weights = somit_ii()
11
12 # SOMIT Part III
13 final_weights = somit_iii(subjective_weights, objective_weights)
```

```
How many criteria (or subcriteria) do you want to compare? 6
You have 6 criteria. Numbered 1 through 6.
->Which criterion (1-6) you think is the median level of importance? 2
--You chose criterion #2 as the median.

Enter importance of each criterion relative to the median.
Compare Criterion #1 with #2: 5
Compare Criterion #3 with #2: 1/3
Compare Criterion #4 with #2:
```

Figure 3. Screenshot of Python code snippet demonstrating the installation and execution of pysomit library.

Additionally, for users who are more comfortable with Excel, we have prepared a spreadsheet template for conducting SOMIT calculations. The template utilizes Excel's built-in Solver tool to solve the constrained optimization problem in SOMIT Part I and employs standard Excel functions for Parts II and III. Both the pysomit Python library and the Excel spreadsheet template are open-source and freely available to all users.

**3.2. Technique for Order of Preference by Similarity to Ideal Solution (TOPSIS)**

The TOPSIS method, an MCDM approach originally proposed by Hwang and Yoon (1981), is designed to evaluate and rank alternatives based on multiple, often conflicting, criteria. TOPSIS is one of the most widely used MCDM methods (Pandey et al., 2023). For instance, Guo and Zhao (2015) employed TOPSIS to rank and select the optimal electric vehicle charging station considering environmental, economic, and social criteria. Elkadeem et al. (2021) used TOPSIS to evaluate and choose the most sustainable energy access option by comparing multiple electrification alternatives, taking 12 performance criteria into account. Zhao et al. (2025) utilized the TOPSIS method is used to select the optimal load regulation for a renewable energy power generation system considering energy efficiency, exergy efficiency, economy, and environmental performance.

In essence, the TOPSIS method is based on the core principle that the most preferred alternative is the one closest to the positive ideal solution (PIS) and farthest from the negative ideal solution (NIS) (Nabavi et al., 2024). The PIS represents the optimal value for each criterion, that is, the highest value for criteria to be maximized and the lowest for those to be minimized. Conversely, the NIS comprises the least desirable values (the lowest for maximization criteria and the highest for minimization criteria). A simple illustration of the five alternatives considering two criteria (one for maximization and one for minimization), along with the ideal solutions, NIS and PIS, is shown in Figure 4.

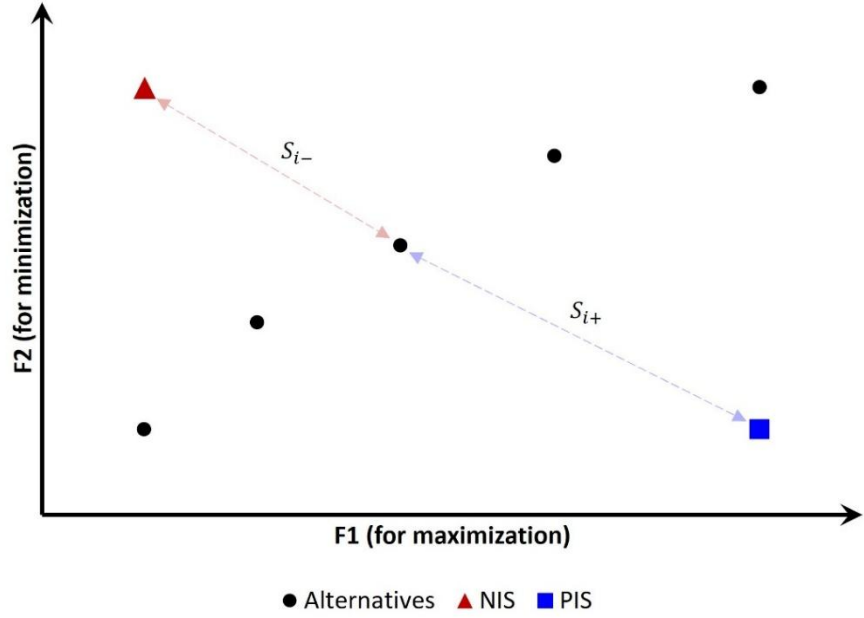

Figure 4. A simple illustration of five alternatives, considering two criteria (F1 for maximization and F2 for minimization), along with the ideal solutions, NIS and PIS.

The following outlines the key steps of the TOPSIS method; refer to Hwang and Yoon (1981) for the original formulation and discussion.

**Step 1:** Normalize the original ACM with $m$ rows (i.e., alternatives) and $n$ columns (i.e., criteria) using the Vector normalization method, as in Hwang and Yoon (1981). It is worth mentioning the row and column indices (i.e., $i$ and $j$) satisfy that $i \in \{1,2,\ldots,m\}$ and $j \in \{1,2,\ldots,n\}$ in all the steps of TOPSIS.

$$F_{ij} = \frac{f_{ij}}{\sqrt{\sum_{k=1}^{m} f_{kj}^2}} \quad (8)$$

where $f_{ij}$ and $f_{kj}$ are the values of the $i^{th}$ and $k^{th}$ alternative under the $j^{th}$ criterion in the original ACM; $F_{ij}$ is the corresponding normalized value in the normalized ACM.

**Step 2:** Construct the weighted normalized ACM by multiplying the values of each criterion with its weight, $w_j$. In this study, the weights are derived using the SOMIT method.

$$v_{ij} = F_{ij} \times w_j \quad (9)$$

**Step 3.** From the weighted normalized ACM, determine the PIS, $A^+$ and NIS, $A^-$, as follows:

$$\begin{aligned} A^+ &= \left\{ \left( \max_i (v_{ij}) \,\middle|\, j \in J_{max} \right), \left( \min_i (v_{ij}) \,\middle|\, j \in J_{min} \right) \right\} \\ &= \{v_1^+, v_2^+, \ldots, v_j^+, \ldots, v_n^+\} \end{aligned} \quad (10)$$

$$A^- = \left\{\left(\min_i(v_{ij}) \middle| j \in J_{max}\right), \left(\max_i(v_{ij}) \middle| j \in J_{min}\right)\right\} \\ = \{v_1^-, v_2^-, \ldots, v_j^-, \ldots, v_n^-\} \tag{11}$$

Where $J_{max}$ is the set of maximization criteria and $J_{min}$ is the set of minimization criteria.

**Step 4.** Compute the Euclidean distances of each alternative to both PIS and NIS, respectively.

$$\text{Distance to PIS}: S_{i+} = \sqrt{\sum_{j=1}^{n}(v_{ij} - v_j^+)^2} \tag{12}$$

$$\text{Distance to NIS}: S_{i-} = \sqrt{\sum_{j=1}^{n}(v_{ij} - v_j^-)^2} \tag{13}$$

**Step 5.** Compute the performance score ($P_i$) of each alternative as follows:

$$P_i = \frac{S_{i-}}{S_{i-} + S_{i+}} \tag{14}$$

Finally, the alternative with the largest $P_i$ is top-ranked and recommended to the decision-maker.

To streamline the implementation of the TOPSIS method and eliminate the need for manual calculations, Pereira et al. (2024) recently developed the pyDecision Python library for multi-criteria decision-making (MCDM). This comprehensive package integrates a wide range of well-established MCDM methods. An example illustrating the use of the pyDecision library will be presented in the next section. In addition, in our previous studies (Wang & Rangaiah, 2017; Wang et al., 2023), we also developed an Excel-based MCDM tool, EMCDM, designed for similar purposes and tailored for decision-makers who are more comfortable working in Excel. It has been shared with a few hundred researchers and practitioners over the years.

## 4. Applications

### 4.1. Renewable Energy Evaluation in India

India, the third-largest energy consumer in the world, is experiencing an ever-growing demand for power due to its expanding population, rapid urbanization, and economic development. Despite having vast conventional energy reserves like coal, oil, and biomass, the country faces a significant energy shortfall, with fossil fuels still meeting about 80% of its demand (International Energy Agency, 2021). This heavy reliance on non-renewable sources has led to concerns over sustainability and energy security. In response, India is making a substantial growth rate in renewable energy. However, each renewable source presents unique financial, technical, environmental, and social challenges, making it difficult to determine the most viable renewable energy solution. These complexities underscore the need for a robust MCDM workflow to guide investment and policy in India's renewable energy sector. In the context of India's renewable energy development, the case study in this section focuses on four key energy sources: solar, wind, hydro, and biomass. These were selected due to their maturity, widespread application, and significant contribution to the country's energy mix. The data of this case study is taken from a recent work by Husain et al. (2024).

**Solar Energy**: Utilizing PV technology, solar energy converts sunlight directly into electricity. India receives an estimated 5,000 trillion kWh of solar energy annually, with most regions experiencing solar irradiance levels ranging from 4 to 7 kWh/m²/day. Owing to its high solar exposure, the country possesses significant solar energy potential, estimated at approximately 748 GW (Husain et al., 2024). The country had installed approximately 60.8 GW of solar capacity, making it one of the fastest-growing solar markets globally.

**Wind Energy**: Wind energy harnesses the kinetic energy of air movement through wind turbines. The country has an installed capacity of 41.6 GW, ranking fourth globally. Evaluations estimate India's wind potential at 302 GW at 100 meters and 695.5 GW at 120 meters elevation (Husain et al., 2024). While current deployment is primarily onshore (which requires valuable land resources that are increasingly becoming a significant constraint), the country also has substantial offshore potential yet to be tapped.

**Hydropower**: This source generates electricity by channeling flowing water through turbines. It contributes roughly 11.5% to India's total electricity generation. With a combined large and small hydro capacity exceeding 50 GW, India is one of the top hydropower producers globally, though environmental and social concerns pose challenges to its expansion.

**Biomass Energy**: Derived from organic sources like agricultural residues and animal waste, biomass energy includes both traditional methods (e.g., combustion of wood, dung, and charcoal) and modern techniques (e.g., biogas and liquid biofuels). India's biomass power capacity stands at 10.2 GW, with significant potential due to the availability of over 750 million metric tonnes of biomass annually (Husain et al., 2024). Its decentralized nature and rural job creation potential make it especially relevant for India's energy landscape.

To evaluate renewable energy sources effectively, this case study considers four principal criterion groups, namely, Financial, Technical, Environmental, and Social, with each encompassing two or more criteria.

**Financial**: This group includes factors directly related to economic investment and operation. The three key criteria are as follows. Total Installed Cost ($/kW): The capital required to purchase and install the energy system; Operation & Maintenance (O&M) Cost ($/kW/y): Ongoing expenses for scheduled and unscheduled maintenance, staffing, and asset management. LCOE, i.e., Levelized Cost of Electricity ($/kWh): Represents the average cost per unit of electricity generated over the system's lifetime, accounting for capital and operating costs.

**Technical**: These assess the functional and performance capabilities of each energy technology. The three key criteria are as follows. Efficiency (%): The percentage of input energy that is converted into usable electricity, a critical factor in system planning. Capacity Factor (%): The actual energy output as a percentage of the system's maximum potential output over time. Technical Maturity (unitless): A qualitative measure rated on a scale from 1 to 5, where 1 indicates the lowest and 5 the highest level of maturity. This rating, provided by domain experts, reflects how well-developed and proven the energy technology is in the current market.

**Environmental**: This group evaluates the ecological impact of each energy source. The two key criteria are as follows. GHG Emissions (g$CO_2$/kWh): The amount of GHG emitted across the full life cycle of the technology. Studies show that the production, installation, and use of renewable energy sources involve a hidden carbon footprint (de Souza Grilo et al., 2018). Land Requirement ($m^2$/kW): The amount of land area needed for deployment, which is a key concern in densely populated regions like India.

**Social**: This group relates to the societal effects and acceptance of renewable energy systems: Job Creation (Job-years/GWh): The ability of a technology to generate employment during both construction and operational phases. Social Acceptance: Assessed through survey responses using a 1–5 qualitative scale, where 1 indicates strong public opposition and 5 indicates strong public approval. This reflects the general level of public opinion and support for the technology.

Moreover, as shown in Table 2, Efficiency, Capacity Factor, Technical Maturity, Job Creation, and Social Acceptance are beneficial criteria, meaning they are for maximization and a higher value indicates a more favorable outcome. In contrast, Total Installed Cost, O&M cost, LCOE, GHG Emissions, and Land Requirement are cost-related criteria, meaning they are for minimization, where lower values are preferred.

Table 2. Ten criteria considered in the case study of renewable energy selection in India.

| Criterion Group | Criterion Name (Unit) | Criterion Code | Type |
|---|---|---|---|
| Financial | Total Installed Cost ($/kW) | C1 | Cost |
|  | O&M Cost ($/kW/y) | C2 | Cost |
|  | LCOE ($/kWh) | C3 | Cost |
| Technical | Efficiency (%) | C4 | Benefit |
|  | Capacity Factor (%) | C5 | Benefit |
|  | Technical Maturity | C6 | Benefit |
| Environmental | GHG Emissions (gCO$_2$/kWh) | C7 | Cost |
|  | Land Requirement (m$^2$/kW) | C8 | Cost |
| Social | Job Creation (Job-years/GWh) | C9 | Benefit |
|  | Social Acceptance | C10 | Benefit |

Husain et al. (2024) extensively reviewed the literature and gathered data from a wide range of credible sources to construct the ACM presented in Table 3. In total, there are four alternatives and ten criteria. Quantitative data for criteria such as Total Installed Cost, O&M cost, Efficiency, and Land Requirement were extracted from existing academic publications, technical reports, and government databases. In addition to literature-based data, the authors also collected qualitative inputs through direct methods: Technical Maturity was assessed using expert judgment on a 1–5 scale, while Social Acceptance was evaluated via a questionnaire survey administered to the public, also rated on a 5-point scale.

Table 3. Alternatives-criteria matrix (ACM) for renewable energy selection in India.

|  | Financial | | | Technical | | | Environmental | | Social | |
|---|---|---|---|---|---|---|---|---|---|---|
|  | C1 | C2 | C3 | C4 | C5 | C6 | C7 | C8 | C9 | C10 |
| Solar | 596 | 9 | 0.038 | 22 | 19 | 4 | 48 | 12 | 0.87 | 4.58 |
| Wind | 1038 | 28 | 0.04 | 35 | 33 | 4 | 11 | 250 | 0.17 | 4.17 |
| Hydro | 1817 | 45.425 | 0.065 | 76.61 | 57 | 5 | 24 | 500 | 0.27 | 3.56 |
| Biomass | 1154 | 46.16 | 0.057 | 84.33 | 68 | 3 | 230 | 13 | 0.21 | 4 |

After obtaining the ACM, the proposed SOMIT method is first applied to determine the weights of the criteria. Subsequently, the TOPSIS method is employed to evaluate and rank the four alternatives considering all ten criteria, using the weights derived from SOMIT.

Since this case study considers four criterion groups (i.e., Financial, Technical, Environmental, and Social), we first apply the steps of SOMIT Part I to assess and calculate the subjective weights assigned to each group. We use $W_1^s$ (in capital form) to denote the subjective weight of Financial group, $W_2^s$ for Technical group, $W_3^s$ for Environmental group, and $W_4^s$ for Social group. Subsequently, these group-level weights are distributed among the individual criteria within each group.

Out of the four criterion groups, the decision-maker selects the Environmental group as having a median level of importance. Each of the remaining criterion groups is then compared against the Environmental group using the 1–9 scale described in Table 1. Here are the subjective assessments: the decision-maker considers that Financial is in between moderately and strongly more important than Environmental, so $a_{13} = 4$; Technical is moderately more important than Environmental, so $a_{23} = 3$; Social is considered to be between equally and moderately less important than Environmental, so $a_{43} = 1/2$. Among the comparisons against Environmental, Financial has the highest comparison value, while Social has the lowest. The decision-maker then directly compares Financial with Social, and judges that Financial is between equally and moderately more important than Social, resulting in $a_{14} = 2$.

With the four comparison values (i.e., $a_{13} = 4$, $a_{23} = 3$, $a_{43} = 1/2$, $a_{14} = 2$), the optimization problem with an equality constraint can be formulated as follows:

$$\begin{aligned} \text{Min} \quad z &= (a_{13}W_3^s - W_1^s)^2 + (a_{23}W_3^s - W_2^s)^2 + (a_{43}W_3^s - W_4^s)^2 + (a_{14}W_4^s - W_1^s)^2 \\ &= (4W_3^s - W_1^s)^2 + (3W_3^s - W_2^s)^2 + \left(\frac{1}{2}W_3^s - W_4^s\right)^2 + (2W_4^s - W_1^s)^2 \end{aligned} \tag{15}$$

$$\begin{aligned} \text{s.t.} \quad & W_1^s + W_2^s + W_3^s + W_4^s - 1 = 0 \\ & 0 \leq W_1^s \leq 1 \\ & 0 \leq W_2^s \leq 1 \\ & 0 \leq W_3^s \leq 1 \\ & 0 \leq W_4^s \leq 1 \end{aligned} \tag{16}$$

Next, use the method of Lagrange multipliers and find all the first-order partial derivatives:

$$\begin{aligned} L &= z - \alpha g \\ &= (4W_3^s - W_1^s)^2 + (3W_3^s - W_2^s)^2 + \left(\frac{1}{2}W_3^s - W_4^s\right)^2 + (2W_4^s - W_1^s)^2 \\ &\quad - \alpha(W_1^s + W_2^s + W_3^s + W_4^s - 1) \end{aligned} \tag{17}$$

where $\alpha$ is the Lagrange multiplier and $g$ represents the equality constraint.

$$\frac{\partial L}{\partial W_1^s} = -2(4W_3^s - W_1^s) - 2(2W_4^s - W_1^s) - \alpha = 0 \tag{18}$$

$$\frac{\partial L}{\partial W_2^s} = 2(W_2^s - 3W_3^s) - \alpha = 0 \tag{19}$$

$$\frac{\partial L}{\partial W_3^s} = 8(4W_3^s - W_1^s) + 6(3W_3^s - W_2^s) + \left(\frac{1}{2}W_3^s - W_4^s\right) - \alpha = 0 \tag{20}$$

$$\frac{\partial L}{\partial W_4^s} = 2\left(W_4^s - \frac{1}{2}W_3^s\right) + 4(2W_4^s - W_1^s) - \alpha = 0 \tag{21}$$

$$\frac{\partial L}{\partial \alpha} = W_1^s + W_2^s + W_3^s + W_4^s - 1 = 0 \tag{22}$$

We now have a system of five linear equations (Eqs. 18 to 22) with five unknowns. Manually solving this system yields the following numerical solution: the weight of Financial is $W_1^s = 0.3900$; the weight of Technical is $W_2^s = 0.3343$; the weight of Environmental is $W_3^s = 0.1056$; and the weight of Social is $W_4^s = 0.1701$.

Alternatively, this optimization problem can be solved directly using our developed pysomit library. The corresponding simple Python code, interactive prompts, and final results are shown in Figure 5. As illustrated, the pysomit library produces the same set of subjective weights as those obtained through manual calculation for the four criterion groups.

```python
# Uncomment next line if you haven't installed the pysomit library
# !pip install pysomit

from pysomit import somit_i, somit_ii, somit_iii

# SOMIT Part I
W_s_group = somit_i()
```

```
How many criteria (or subcriteria) do you want to compare? 4
You have 4 criteria. Numbered 1 through 4.
->Which criterion (1-4) you think is the median level of importance? 3
--You chose criterion #3 as the median.

Enter importance of each criterion relative to the median.
Compare Criterion #1 with #3: 4
Compare Criterion #2 with #3: 3
Compare Criterion #4 with #3: 1/2

Highest: criterion #1 (4.0); Lowest: criterion #4 (0.5)
Compare Criterion #1 with #4: 2
 a_hl = 2.0

Optimal subjective weights (w^s):
  w^s_1 = 0.3900
  w^s_2 = 0.3343
  w^s_3 = 0.1056
  w^s_4 = 0.1701
Objective value z = 0.017595
```

Figure 5. Screenshot of Python code snippet using the pysomit library for solving the formulated optimization problem in SOMIT Part I.

Subsequently, the Financial group comprises three criteria: Total Installed Cost, O&M Cost, and LCOE. These internal criteria are compared pairwise to determine their respective shares (fractions) of the total subjective weight allocated to the Financial group. The decision-maker first selects O&M Cost as the median importance criterion within this group. As SOMIT flexibly supports the use of decimal values for finer granularity (as shown in the last row of Table 1), the decision-maker assesses that: Total Installed Cost is considered 1.5 times as important as O&M Cost; LCOE is judged to be 0.8 times as important as O&M Cost. Given that Total Installed Cost has the highest comparison value and LCOE the lowest, the decision-maker directly compares these two criteria and assesses that Total Installed Cost is 1.3 times as important as LCOE. Using these values, a constrained optimization problem (similar to Eqs. 15 and 16) is formulated and solved using the method of Lagrange multipliers or other suitable computational tools (e.g., pysomit). The resulting relative shares within the Financial group are Total Installed Cost: 0.4135, O&M Cost: 0.2966, LCOE: 0.2899. Accordingly, their subjective weights are computed by multiplying their respective shares by the allocated group weight $W_1^s = 0.3900$. Thus, the resulting subjective weights are as follows, Total Installed Cost: $w_1^s = 0.3900 \times 0.4135 = 0.1613$, O&M Cost: $w_2^s = 0.3900 \times 0.2966 = 0.1157$, LCOE: $w_3^s = 0.3900 \times 0.2899 = 0.1131$.

For the remaining three criterion groups (Technical, Environmental, and Social), their subjective group weights are evenly distributed among the criteria within each group. This is equivalent to assigning equal importance (a value of 1) when comparing the criteria (e.g., when using pysomit). The complete set of subjective weights for all criteria is presented in the second column of Table 4.

Table 4. Criteria weights determined using SOMIT for renewable energy selection in India.

| Criterion Code | Subjective Weight $w_j^s$ | Objective Weight $w_j^o$ | Final Weight $w_j$ |
| --- | --- | --- | --- |
| C1 | 0.1613 | 0.0830 | 0.1335 |
| C2 | 0.1157 | 0.1113 | 0.1285 |
| C3 | 0.1131 | 0.1235 | 0.1393 |
| C4 | 0.1114 | 0.1263 | 0.1405 |
| C5 | 0.1114 | 0.1129 | 0.1255 |
| C6 | 0.1114 | 0.0758 | 0.0843 |
| C7 | 0.0528 | 0.0841 | 0.0443 |
| C8 | 0.0528 | 0.1126 | 0.0593 |
| C9 | 0.0850 | 0.0823 | 0.0698 |
| C10 | 0.0850 | 0.0884 | 0.0750 |

Next, following the steps of SOMIT Part II, the objective weights are derived from the ACM directly. For brevity, the calculation details are provided in Table 5; while the full set of objective weights for all criteria is summarized in the third column of Table 4, after applying Eq. 6. Proceeding to SOMIT Part III, using the Eq. 7, the subjective and objective weights are then integrated to obtain the final weights, which are presented in the last column of Table 4. These weights can also be visualized as bar charts in Figure 6.

Table 5. Calculation details of SOMIT Part II for renewable energy selection in India.

| | Financial | | | Technical | | | Environmental | | Social | |
|---|---|---|---|---|---|---|---|---|---|---|
| | C1 | C2 | C3 | C4 | C5 | C6 | C7 | C8 | C9 | C10 |
| Normalized ACM (using Max-Min Normalization) | | | | | | | | | | |
| Solar | 0 | 0 | 0 | 0 | 0 | 0.5 | 0.1689 | 0 | 1 | 1 |
| Wind | 0.3620 | 0.5113 | 0.0741 | 0.2086 | 0.2857 | 0.5 | 0 | 0.4877 | 0 | 0.5980 |
| Hydro | 1 | 0.9802 | 1 | 0.8761 | 0.7755 | 1 | 0.0594 | 1 | 0.1429 | 0 |
| Biomass | 0.4570 | 1 | 0.7037 | 1 | 1 | 0 | 1 | 0.0020 | 0.0571 | 0.4314 |
| Median | 0.4095 | 0.7458 | 0.3889 | 0.5424 | 0.5306 | 0.5 | 0.1142 | 0.2449 | 0.1 | 0.5147 |
| Absolute Deviation from Median (ADM) | | | | | | | | | | |
| Solar | 0.4095 | 0.7458 | 0.3889 | 0.5424 | 0.5306 | 0.0000 | 0.0548 | 0.2449 | 0.9000 | 0.4853 |
| Wind | 0.0475 | 0.2345 | 0.3148 | 0.3338 | 0.2449 | 0.0000 | 0.1142 | 0.2428 | 0.1000 | 0.0833 |
| Hydro | 0.5905 | 0.2345 | 0.6111 | 0.3338 | 0.2449 | 0.5000 | 0.0548 | 0.7551 | 0.0429 | 0.5147 |
| Biomass | 0.0475 | 0.2542 | 0.3148 | 0.4576 | 0.4694 | 0.5000 | 0.8858 | 0.2428 | 0.0429 | 0.0833 |
| AADM ($r_j$) | 0.2738 | 0.3672 | 0.4074 | 0.4169 | 0.3724 | 0.2500 | 0.2774 | 0.3714 | 0.2714 | 0.2917 |

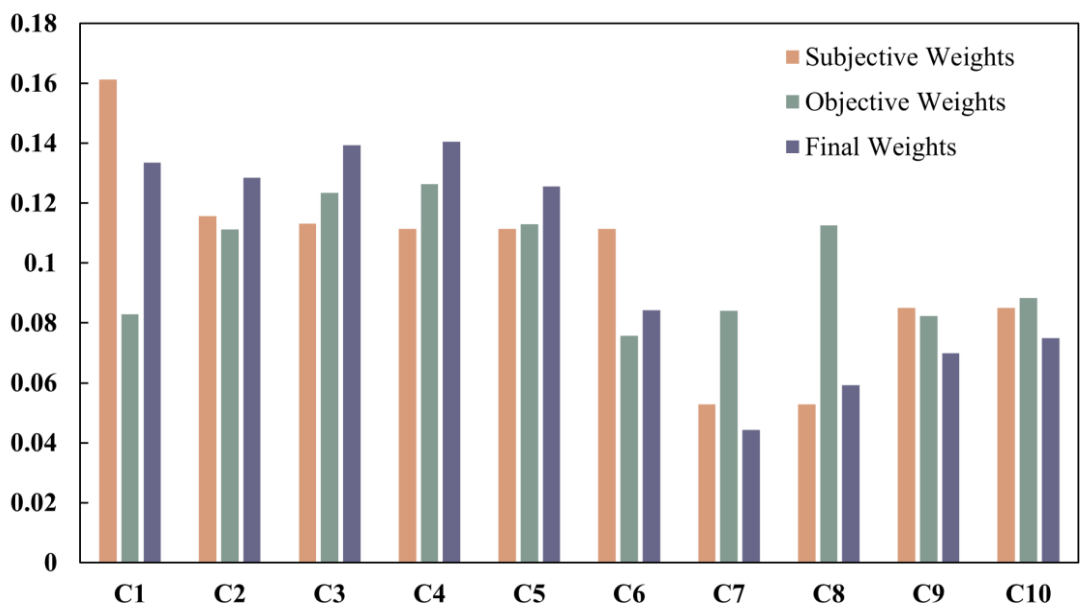

Figure 6. Subjective, objective, and final weights determined by SOMIT for the ten criteria of renewable energy selection in India.

The complete code using our developed pysomit library to solve all three parts of the SOMIT method in a single run (eliminating the need for manual calculations) is shown in Table 6. The output subjective, objective, and final combined weights are displayed as comments in the last few lines of Table 6, following the hash symbol (#). These generated weights are consistent with those obtained through manual calculation.

Table 6. Complete code using pysomit library to find the weights of criteria for renewable energy selection in India.

```python
# Uncomment next line if you haven't installed the pysomit library
# !pip install pysomit

from pysomit import somit_i, somit_ii, somit_iii
import numpy as np

# SOMIT Part I
W_s_group = somit_i()
F_frac = somit_i()
T_frac = somit_i()
E_frac = somit_i()
S_frac = somit_i()

groups = [F_frac, T_frac, E_frac, S_frac]
w_s = [W * frac for W, gp in zip(W_s_group, groups) for frac in gp]
w_s = np.array(w_s, dtype=float)
print("\nSubjective Weights:",  np.round(w_s,4))

# SOMIT Part II
ACM = np.array([
    [596,   9,      0.038, 22,    19, 4,    48,  12,  0.87, 4.58],
    [1038,  28,     0.04,  35,    33, 4,    11,  250, 0.17, 4.17],
    [1817,  45.425, 0.065, 76.61, 57, 5,    24,  500, 0.27, 3.56],
    [1154,  46.16,  0.057, 84.33, 68, 3,    230, 13,  0.21, 4.00]
])

w_o = somit_ii(ACM)
print("\nObjective Weights:", np.round(w_o,4))

# SOMIT Part III
final_weights = somit_iii(w_s, w_o)
print("\nFinal Weights:", np.round(final_weights,4))

# -------------- Output Weights ----------------
# Subjective Weights: [0.1613 0.1157 0.1131 0.1114 0.1114
#                     0.1114 0.0528 0.0528 0.085  0.085 ]
# Objective Weights: [0.083  0.1113 0.1235 0.1263 0.1129
```

```
#                       0.0758 0.0841 0.1126 0.0823 0.0884]
# Final Weights: [0.1335 0.1285 0.1393 0.1405 0.1255
#                 0.0843 0.0443 0.0593 0.0698 0.075 ]
# --------------------------------------------
```

Next, using the final weights determined by SOMIT, the TOPSIS method is applied to rank the four alternatives across all ten criteria. Following the procedure outlined in Section 2.2, the weighted normalized ACM, along with PIS and NIS, are presented in Table 7. Note that vector normalization is used here, as in Hwang and Yoon (1981). It is worth reiterating that Efficiency, Capacity Factor, Technical Maturity, Job Creation, and Social Acceptance are maximization criteria, whereas Total Installed Cost, O&M cost, LCOE, GHG Emissions, and Land Requirement are minimization criteria. This distinction should be carefully taken into consideration when identifying the PIS and NIS in Step 3 of the TOPSIS method.

Table 7. Weighted normalized ACM and ideal solutions for ranking alternatives using the TOPSIS method.

| | Weighted Normalized ACM (using Vector Normalization) | | | | | | | | | |
| --- | --- | --- | --- | --- | --- | --- | --- | --- | --- | --- |
| | Financial | | | Technical | | | Environmental | | Social | |
| | C1 | C2 | C3 | C4 | C5 | C6 | C7 | C8 | C9 | C10 |
| Solar | 0.0323 | 0.0163 | 0.0516 | 0.0255 | 0.0247 | 0.0415 | 0.0090 | 0.0013 | 0.0639 | 0.0420 |
| Wind | 0.0563 | 0.0506 | 0.0543 | 0.0406 | 0.0429 | 0.0415 | 0.0021 | 0.0265 | 0.0125 | 0.0382 |
| Hydro | 0.0985 | 0.0821 | 0.0883 | 0.0888 | 0.0741 | 0.0519 | 0.0045 | 0.0530 | 0.0198 | 0.0326 |
| Biomass | 0.0626 | 0.0834 | 0.0774 | 0.0978 | 0.0884 | 0.0311 | 0.0431 | 0.0014 | 0.0154 | 0.0366 |
| | Ideal Solutions | | | | | | | | | |
| PIS | 0.0323 | 0.0163 | 0.0516 | 0.0978 | 0.0884 | 0.0519 | 0.0021 | 0.0013 | 0.0639 | 0.0420 |
| NIS | 0.0985 | 0.0834 | 0.0883 | 0.0255 | 0.0247 | 0.0311 | 0.0431 | 0.0530 | 0.0125 | 0.0326 |

Subsequently, the Euclidean distances of each alternative to the PIS and NIS are calculated and presented in the second and third rows of Table 8, respectively. Based on these values and using the Eq. 14, the performance score of each alternative is computed and shown in the last column of Table 8.

Table 8. Distances to PIS and NIS, as well as performance scores

| Alternative | $S_{i+}$ | $S_{i-}$ | $P_i$ |
| --- | --- | --- | --- |
| Solar | 0.0971 | 0.1300 | 0.5725 |
| Wind | 0.1025 | 0.0842 | 0.4511 |
| Hydro | 0.1227 | 0.0918 | 0.4279 |
| Biomass | 0.1029 | 0.1157 | 0.5293 |

As can be seen, the TOPSIS results yield the following ranking: Solar ≻ Biomass ≻ Wind ≻ Hydro. This outcome recommends Solar as the top-ranked renewable energy alternative in India, based on all ten criteria and using the SOMIT weights that integrate both subjective assessments and objective data-driven insights. This outcome is consistent with India's broader renewable energy policy landscape, where solar power has been identified as the backbone of the clean energy transition under the National Solar Mission and subsequent state-level initiatives (Umar & Yadaw, 2024). The SOMIT-based prioritization of Solar not only reflects its lower cost (in installation, O&M, and LCOE), but also its comparatively modest land requirements, which is a significant advantage in a densely populated country. The recent review article by Umar & Yadaw (2024) recommends that India should continue strengthening policy frameworks, encouraging private sector involvement, and advancing research and development to fully realize its solar energy potential. However, at the same time, the SOMIT-based results underscore an important trade-off. While solar excels in cost and land use performance, it comes with relatively higher GHG emissions compared with wind and hydro. These emissions are primarily associated with the upstream manufacturing processes, particularly the energy-intensive production of PV modules, rather than operational phases. By capturing this nuance, SOMIT ensures that environmental sustainability is assessed holistically, considering life cycle impacts alongside operational benefits. This provides policymakers with a clearer picture of the trade-offs involved, ensuring that the rapid deployment of solar is balanced with strategies for reducing embedded carbon in manufacturing and supply chains. Moreover, the weighting stage captures the importance of social acceptance, a factor that has influenced the pace of renewables adoption in India (Kar et al., 2024). For example, public opposition to land acquisition and biomass collection logistics has often slowed project implementation despite technical potential. The SOMIT framework also makes these sensitivities visible, so that renewable energies aim to be socially acceptable/sustainable.

Likewise, as presented in Table 9, by applying the TOPSIS function from the pyDecision library coupling with the final weights determined by SOMIT, the same performance scores as those shown in Table 8 are obtained. Employing these reliable tools can avoid the need for manual computation.

As shown, although the final hybrid weights are used in TOPSIS in this case, the modular design of SOMIT allows for flexibility: either the subjective or objective weights can be used independently in the MCDM process, depending on the preference of the decision-maker. Some decision-makers may prioritize their own subjective assessments over weights derived from the ACM. Others may prefer to rely entirely on the objective weights calculated from the ACM in order to minimize subjectivity introduced at this stage. In addition, SOMIT weights can be seamlessly integrated into any MCDM method that requires criteria weights as part of its algorithm. From a managerial perspective, SOMIT provides decision-makers with a structured approach to balance technical, financial, environmental, and social trade-offs in renewable energy project evaluation. In the present case study of India, where renewable energy expansion must address affordability and large-scale deployment, SOMIT enables policymakers and energy managers to prioritize technologies not only based on cost and efficiency but also on broader socio-economic impacts such as job creation and social acceptance. This helps managers align investment decisions with both national renewable targets and developmental objectives.

Table 9. Python code of using the TOPSIS method from pyDecision library for ranking renewable energy sources in India.

```python
# Uncomment next line if you haven't installed the pyDecision library
# !pip install pyDecision

import numpy as np
from pyDecision.algorithm import topsis_method

# Criterion type: 'max' for benefit;  'min' for cost
types = ['min', 'min', 'min', 'max', 'max', 'max',
'min','min','max','max']

# Rows: Alternatives; Columns: Criteria
ACM = np.array([
    [596,   9,      0.038, 22,    19, 4,   48,  12,  0.87, 4.58],
    [1038,  28,     0.04,  35,    33, 4,   11,  250, 0.17, 4.17],
    [1817,  45.425, 0.065, 76.61, 57, 5,   24,  500, 0.27, 3.56],
    [1154,  46.16,  0.057, 84.33, 68, 3,   230, 13,  0.21, 4.00]
])

# Here, take the hybrid final_weights found by SOMIT
weights = np.round(final_weights,4)

# Calculation using TOPSIS method
performance_scores  = topsis_method(ACM, weights, types,
                                    graph=False, verbose=False)
print("Performance Scores:", np.round(performance_scores,4))

# -------------- Output -----------------
# Performance Scores: [0.5725 0.4511 0.4279 0.5293]
# ---------------------------------------
```

## 4.2. Evaluation of Renewable Energy Sources in Saudi Arabia

To further demonstrate the practicality of the proposed SOMIT method in supporting MCDM, this section briefly presents another case study focused on evaluating renewable energy resources in Saudi Arabia, adopted from Yazdani et al. (2020). Five renewable energy sources, namely, Solar PV (photovoltaic), Solar Thermal, Wind, Geothermal, and Biomass are evaluated and ranked based on a set of criteria that capture technical, financial, environmental, and social dimensions. Eight key criteria (C1 to C8 as follows) were selected from the existing literature to guide this decision analysis. The complete ACM data, taken from Yazdani et al. (2020), is presented in Table 10.

- C1: Capital Cost (USD/MW) accounts for all expenses related to establishing a power plant, including land acquisition, equipment procurement, labor, installation, and infrastructure development. This is cost criterion.
- C2: O&M Cost (USD/KW-year) encompasses recurring expenditures such as staffing, routine servicing, and costs associated with energy generation. This is cost criterion.
- C3: Energy Cost (USD/kWh) reflects the projected cost of electricity production over the plant's operational lifetime. This is cost criterion.
- C4: Job Creation (Total job-year/GWh) measures the annual employment potential associated with each renewable energy system. This is a benefit criterion.
- C5: Land Use ($m^2$/GWh) quantifies the physical space required to construct and operate the energy facility. This is cost criterion.
- C6: GHG emissions ($tCO_2$ equivalent/MWh) represent the amount of harmful gases released by the plant during operation, providing an environmental impact indicator. This is cost criterion.
- C7: Efficiency (%) indicates how effectively the energy source is converted into usable electricity, expressed as a percentage of energy output relative to input. This is a benefit criterion.
- C8: resource availability (kwh/$m^2$/year) assesses the extent to which renewable resources are locally accessible and sufficient to support consistent electricity generation. This is a benefit criterion.

Table 10. ACM for evaluation of renewable energy sources in Saudi Arabia.

|  | C1 | C2 | C3 | C4 | C5 | C6 | C7 | C8 |
|---|---|---|---|---|---|---|---|---|
| Solar PV | 3873 | 39.55 | 0.27 | 0.87 | 150 | 0.07 | 12 | 2130 |
| Solar Thermal | 5067 | 67.26 | 0.23 | 0.23 | 40 | 0.02 | 21 | 2200 |
| Wind | 2213 | 24.69 | 0.08 | 0.17 | 200 | 0.04 | 35 | 570 |
| Geothermal | 6243 | 132 | 0.07 | 0.25 | 100 | 0.04 | 16 | 100 |
| Biomass | 8312 | 460.47 | 0.05 | 0.21 | 25 | 0.1 | 25 | 200 |

In this instance, the SOMIT method is firstly applied to directly assess the eight criteria (unlike the previous case study, which first assesses criterion groups before examining individual criteria within each group). Here, Criterion C3 is taken as the median for pairwise comparisons. The subjective assessments were as follows: C1 vs. C3 = 3; C2 vs. C3 = 2; C4 vs. C3 = 1; C5 vs. C3 = 1/4; C6 vs. C3 = 1/2; C7 vs. C3 = 1; and C8 vs. C3 = 3. Based on these assessments, C1 gets the highest comparison value, while C5 has the lowest. An additional comparison of C1 vs. C5 is conducted, assigning a value of 3. The subjective, objective, and final combined weights derived using the SOMIT method are summarized in Table 11. For brevity, intermediate calculations are not shown but can be easily reproduced using the pysomit library or Excel Solver.

Table 11. Weights of all criteria determined using SOMIT for evaluation of renewable energy sources in Saudi Arabia.

| Criterion Code | Subjective Weight $w_j^s$ | Objective Weight $w_j^o$ | Final Weight $w_j$ |
|---|---|---|---|
| C1 | 0.2289 | 0.1187 | 0.2111 |
| C2 | 0.1642 | 0.1036 | 0.1322 |
| C3 | 0.0806 | 0.1476 | 0.0924 |
| C4 | 0.0836 | 0.0904 | 0.0587 |
| C5 | 0.0710 | 0.1392 | 0.0768 |
| C6 | 0.0433 | 0.1175 | 0.0396 |
| C7 | 0.0836 | 0.1189 | 0.0773 |
| C8 | 0.2448 | 0.1640 | 0.3120 |

It can be observed that only eight pairwise comparisons are required for the eight criteria (i.e., a linear relationship with respect to $n$ criteria). In contrast, the conventional AHP method would require 28 comparisons, following the quadratic expression of $n(n-1)/2$. Next, using the final weights determined by SOMIT, the TOPSIS method is applied to evaluate the alternative energy sources. The resulting performance scores are 0.7562 for Solar PV, 0.7556 for Solar Thermal, 0.4997 for Wind, 0.3312 for Geothermal, and 0.2281 for Biomass. Thus, the ranking is: Solar PV ≻ Solar Thermal ≻ Wind ≻ Geothermal ≻ Biomass. These results indicate that Solar PV is the most preferred renewable energy source in Saudi Arabia, considering all eight criteria simultaneously. They further align with Saudi Arabia's Vision 2030 program, which designates solar energy as the cornerstone of its renewable portfolio, targeting 58.7 GW of renewable capacity, of which 40 GW (i.e., 68.1%) is expected from solar (Ali et al., 2021). By confirming the dominance of solar PV and solar thermal within a robust MCDM framework, SOMIT strengthens the analytical basis for ongoing energy policy and decision-making in Saudi Arabia's renewable energy landscape. At the same time, the SOMIT-based results also underscore why alternatives such as geothermal and biomass rank lower in this context (ranked second-to-last and last, respectively). Geothermal resources in Saudi Arabia are geographically limited and not cost-competitive at scale, while biomass faces constraints in terms of resource availability and long-term sustainability in arid regions. These disadvantages, clearly reflected in the weighting and evaluation MCDM process, explain their marginal roles in current policy plans. By bringing such insights to light, the proposed framework not only confirms the suitability of solar power but also illustrates the challenges of diversifying the renewable mix beyond solar and wind under local resource conditions. Notably, solar and wind are projected to account for 95.4% of Saudi Arabia's targeted renewable capacity. Here, the managerial relevance of SOMIT lies in its ability to guide decision-making under the country's ambitious energy transition and diversification agenda. SOMIT supports energy planners in weighing the environmental benefits of renewables against the need for economic competitiveness and energy security, thus offering a transparent framework to justify technology prioritization to stakeholders. Moving forward, managers in utilities and ministries can further apply SOMIT to fine-tune optimal technology portfolios (e.g., adjusting the distribution of renewables as technologies evolve and conditions change) that balance the overall carbon neutrality goals with economic diversification strategies under Vision 2030 program.

## 5. Discussion and Comparisons

For comparative purposes, the widely used AHP method is also applied to derive weights of the ten criteria in the India's renewable energy evaluation case study. As in Section 4.1, we first determine the weights of the four criterion groups (Financial, Technical, Environmental, and Social), and then further split the group weights into their constituent criteria, in a manner consistent with SOMIT Part I. Table 12 presents the complete pairwise comparison matrix of the four criterion groups used in AHP method. As shown, the comparison values in the Environmental column are retained from Section 4.1, where the Environmental group was considered to represent a median level of importance. Similarly, the comparison value between the Financial and Social groups is also preserved in the AHP matrix.

Table 12. Full pairwise comparison matrix for the four criteria groups of renewable energy selection in India, using AHP method.

| Criterion Group | Financial | Technical | Environmental | Social |
| --- | --- | --- | --- | --- |
| Financial | 1 | 2 | 4 | 2 |
| Technical | 1/2 | 1 | 3 | 5 |
| Environmental | 1/4 | 1/3 | 1 | 2 |
| Social | 1/2 | 1/5 | 1/2 | 1 |

Solving the AHP problem using the original eigenvalue method with the aid of Python NumPy library, the group weights were obtained as follows: Financial $W_1^s = 0.4257$, Technical $W_2^s = 0.3401$, Environmental $W_3^s = 0.1297$, Social $W_4^s = 0.1045$. Within the Financial group, the same internal comparisons from SOMIT Part I are adopted and solved using the AHP method, outputting relative shares of 0.4108 for Total Installed Cost, 0.3095 for O&M Cost, and 0.2797 for LCOE. Accordingly, their weights are computed by multiplying their respective shares by the Financial group weight. For the remaining three criterion groups (Technical, Environmental, and Social), their respective group weights were distributed equally among the criteria within each group, as in SOMIT Part I. The complete set of AHP weights for all the ten criteria are C1: 0.1749, C2: 0.1318, C3: 0.1191, C4: 0.1134, C5: 0.1134, C6: 0.1134, C7: 0.0649, C8: 0.0649, C9: 0.0523, C10: 0.0523. When compared with the SOMIT subjective weights reported in the second column of Table 4, these AHP weights do not deviate much. This consistency is further supported by their Pearson correlation coefficient of 0.898, which indicates a strong positive correlation between the two sets of subjective weights.

However, as briefly noted in the earlier sections, the conventional AHP method requires *n(n−1)/2* pairwise comparisons to construct the full comparison matrix, where *n* is the number of criteria. This quadratic growth will quickly become burdensome to decision-makers as *n* increases. By contrast, the proposed SOMIT method requires only *n* subjective comparisons in total from decision-makers and is capable of producing comparable weights. Figure 7 clearly illustrates this difference; AHP grows quadratically, and SOMIT increases linearly. Consequently, while largely preserving the quality of the subjective weights, SOMIT substantially reduces the cognitive load on decision-makers, especially in problems involving many criteria (e.g., more than six), which is common in the energy system MCDM applications.

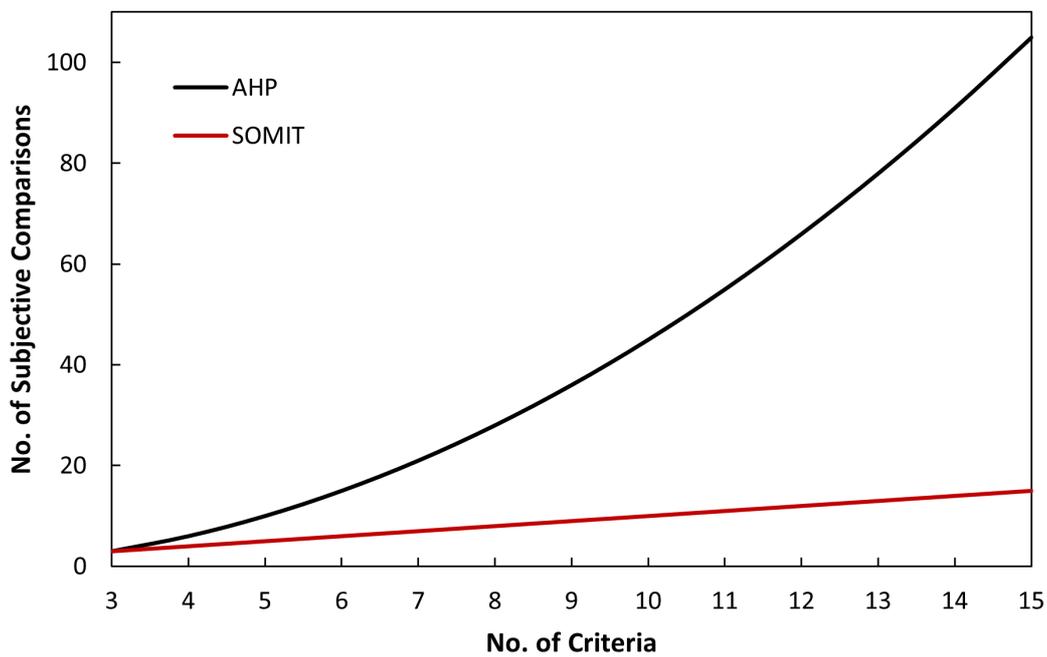

Figure 7. Number of subjective pairwise comparisons required for AHP and SOMIT methods.

Next, we evaluate SOMIT Part II (that derives the objective weights from ACM) and the widely adopted objective weighting method, CRITIC, when confronted with extreme outlier value(s) in ACM. Outliers may be included in the ACM for a variety of reasons, such as measurement errors, data collection or processing issues, and data entry/recording mistakes. For illustration, consider the case study of India's renewable energy evaluation again, and suppose that the $CO_2$ emissions of solar power were incorrectly recorded as 480 $gCO_2$/kWh, instead of the correct value of 48 $gCO_2$/kWh. It is worth noting that an emission of 480 $gCO_2$ per kWh across the full life cycle of solar power is unrealistically high; the average is reported as 49.91 $gCO_2$/kWh (Nugent & Sovacool, 2014). This error thus represents a clear outlier in ACM. We then assess how the presence of this outlier affects the objective criteria weights determined by SOMIT Part II and by CRITIC. As shown in Figure 8, CRITIC is more affected by the extreme value, whereas SOMIT Part II is comparatively less impacted. This difference arises from their underlying principles: SOMIT Part II employs the medians in its calculations, which inherently are more robust toward outliers; while CRITIC relies on arithmetic means for measuring criterion dispersion and pairwise correlations, which are more easily distorted by extreme values. Beyond fewer required subjective comparisons and greater robustness to outliers, an additional advantage of SOMIT is that its Part III always strives to enhance consistency and balance by combining subjective weights from decision-makers with objective data-driven weights from ACM, thus reducing the bias that can result from relying solely on either subjective or objective weighting methods.

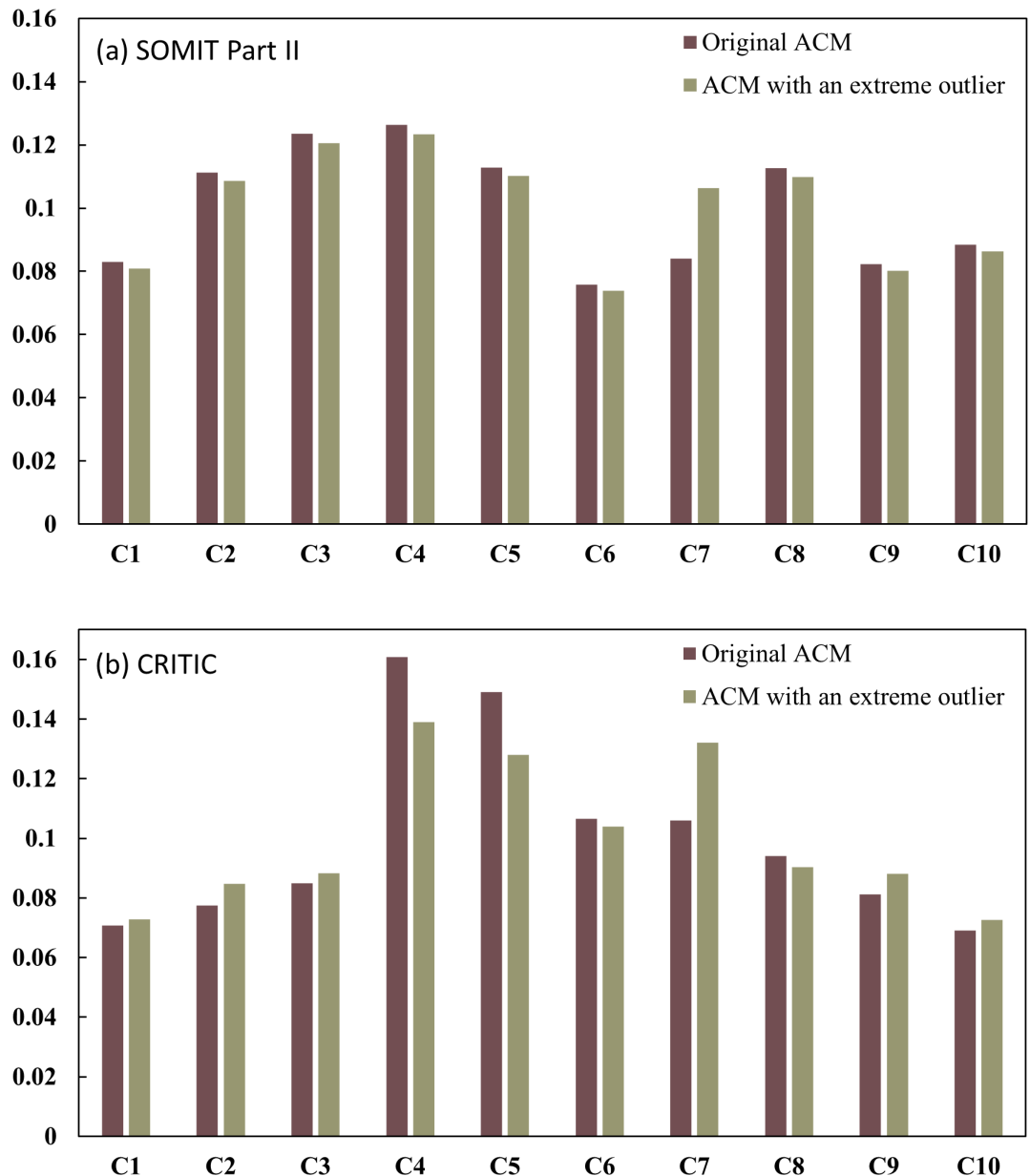

Figure 8. Objective criteria weights determined by (a) SOMIT Part II and (b) CRITIC, based on the original ACM and the ACM with an extreme outlier.

Further, for more sensitivity/robustness analyses, we adopt the comprehensive approach of Nabavi et al. (2023), which evaluates how weighting methods respond to changes or modifications introduced in ACM. The first modification considered is the linear transformation of a criterion, rooted in the principle of independence of the value scale (IVS) for consistency analysis. IVS states that the outcomes of an objective weighting method should remain minimally affected by a linear change in the measurement scale of a criterion. In simple terms, whether a criterion (e.g., social acceptance questionnaire) is gauged on the original 1–5 scale (1, 2, 3, 4, 5) or linearly transformed to a 1–9 scale by a function such as $y = 2x - 1$, the derived outcomes in a rational decision-making process (e.g., weighting step) should be minimally impacted. The second modification is the reciprocal reformulation of a criterion, which may occur during data measurement/collection or optimization model design. For instance, the minimization of GHG emissions ($gCO_2$/kWh) can be equivalently reformulated as the maximization of energy yield per unit of GHG emission (kWh/$gCO_2$) by simply taking the reciprocal, or as the maximization of "a constant minus GHG emissions", where the constant can represent the maximum allowable emission level. In such cases, a rational decision-making framework should also yield outcomes that are minimally impacted. In the India case study on renewable energy evaluation, we introduce multiple modifications to the ACM, as shown in Table 13. First, the Technical Maturity criterion (C6) is linearly transformed from the original 1–5 scale to a 1–9 measurement scale. Second, the GHG Emissions ($gCO_2$/kWh) criterion (C7) is equivalently reformulated as maximizing Energy Yield per GHG Emission (kWh/$gCO_2$) by taking the reciprocal. Third, the Job Creation (Job-years/GWh) criterion, where all the original values are less than 1, is redefined as minimizing "1 minus Job Creation", representing how far each alternative is from the interim goal of achieving 1 Job-year at the current stage of development (with smaller value/distance to the goal indicating better performance).

Table 13. ACM with modifications for sensitivity analysis.

|  | Financial | | | Technical | | | Environmental | | Social | |
| --- | --- | --- | --- | --- | --- | --- | --- | --- | --- | --- |
|  | C1 | C2 | C3 | C4 | C5 | C6 | C7 | C8 | C9 | C10 |
| Solar | 596 | 9 | 0.038 | 22 | 19 | 7 | 1/48 | 12 | 0.13 | 4.58 |
| Wind | 1038 | 28 | 0.04 | 35 | 33 | 7 | 1/11 | 250 | 0.83 | 4.17 |
| Hydro | 1817 | 45.425 | 0.065 | 76.61 | 57 | 9 | 1/24 | 500 | 0.73 | 3.56 |
| Biomass | 1154 | 46.16 | 0.057 | 84.33 | 68 | 5 | 1/230 | 13 | 0.79 | 4 |

Nabavi et al. (2023) described a specific metric, average absolute fractional deviation for weights ($AAFDw$), to quantify changes in weights when the original ACM is modified:

$$AAFDw = \frac{1}{n}\sum_{j=1}^{n}\frac{|WO_j - WM_j|}{\left|\frac{WO_j + WM_j}{2}\right|} \qquad (23)$$

where $n$ is the number of criteria, and $WO_j$ and $WM_j$ are the computed weights for $j^{th}$ criterion in the original and modified ACM, respectively. A minimal sensitivity level is indicated when the $AAFDw$ value approaches zero or remains negligible. Applying this metric, the objective weights obtained by both SOMIT and CRITIC are compared against their respective weights derived from the original ACM. The analysis shows that the $AAFDw$ of SOMIT is 1.9%, while that of CRITIC is 2.5%, demonstrating that SOMIT exhibits higher robustness than the benchmark CRITIC method.

In addition to the methodological comparisons with AHP and CRITIC, the findings of this work carry several managerial and policy implications for renewable energy planning. First, SOMIT provides decision-makers with a structured and transparent approach to balance subjective expert opinions with objective ACM dataset. By requiring fewer subjective comparisons and employing medians to mitigate the influence of outliers, SOMIT enhances both efficiency and robustness. This makes it particularly suitable for large-scale renewable energy assessments, where diverse criteria must be evaluated and data uncertainties are frequent. From a managerial perspective, the case study of India illustrates how SOMIT can highlight not only the cost advantages of solar power but also its trade-offs, such as higher life-cycle GHG emissions relative to other options. These insights equip planners with a more nuanced understanding of deployment challenges, including environmental sustainability and social acceptance, which are critical for achieving the country's ambitious renewable energy targets. In Saudi Arabia, SOMIT supports the strategic prioritization of solar power under Vision 2030, while also clarifying why alternatives such as biomass and geothermal remain less feasible in the local context. For managers and policymakers, this ensures that investment decisions are not only technically sound but also aligned with national priorities and resource realities. More broadly, SOMIT equips energy managers, policymakers, and organizational leaders with a defensible tool for resource allocation, project selection, and policy design across a wide range of energy system applications. Its ability to integrate multiple perspectives, reduce cognitive burden on decision-makers, and yield results that are robust and interpretable provides managers with a practical method for guiding renewable energy strategies. As a result, SOMIT bridges the gap between methodological development and real-world application, contributing directly to evidence-based energy planning and managerial decision-making.

Last but not least, while SOMIT reduces comparison burden, enhances robustness to outliers, and balances subjective and objective inputs, it still relies on the quality and representativeness of the available data and expert judgments. If the ACM of an application is constructed from incomplete, outdated, or biased data, or if the subjective assessments are provided by individuals lacking sufficient domain expertise, then the resulting weights and the subsequent rankings may still be affected. Simply put, although the use of medians mitigates extreme distortions, SOMIT cannot fully eliminate dependence on input quality. Future work could further address this limitation by integrating group decision-making mechanisms and incorporating advanced fuzzy extensions (Baydaş et al., 2025), to better capture uncertainties and imprecisions.

**6. Conclusion**

In conclusion, this study proposed a novel hybrid weighting method, SOMIT, to aid MCDM process. Recognizing the pivotal role of criteria weights, SOMIT addresses limitations in conventional methods by integrating subjective expert judgments with objective data insights in a modular, user-friendly, and computationally efficient framework. The three-part design of SOMIT significantly reduces the cognitive burden on decision-makers by reducing the number of subjective pairwise comparisons to a linear scale, utilizes normalized median-based dispersion to extract meaningful objective weights from the ACM data (that is inherently more robust to outliers), and merges both sources of weights using a multiplicative synthesis to ensure balance and mitigate bias. A dedicated Python library, pysomit, was developed to simplify and streamline the application of the SOMIT method in renewable energy evaluation applications.

SOMIT was applied to two real-world case studies involving the evaluation and selection of renewable energy sources in India and Saudi Arabia. In both cases, SOMIT derived a set of weights for criteria spanning financial, technical, environmental, and social dimensions. These weights were then used in conjunction with TOPSIS method to rank the renewable energy alternatives, considering multiple criteria simultaneously. SOMIT offers a modular and versatile framework that can be readily integrated with a wide array of MCDM methods, hence extending its utility to other sectors and decision problems that require balanced and context-sensitive weighting. From a practical perspective, SOMIT empowers policymakers, planners, and stakeholders with a more transparent and structured approach to deriving criteria weights, an essential step in energy decision-making, particularly in complex and rapidly evolving contexts such as the renewable energy transition. By integrating both expert knowledge and empirical data, SOMIT provides robust, interpretable, and context-sensitive insights that directly support applied energy challenges such as carbon neutrality, technology deployment strategies, and renewable energy investment prioritization. Consequently, the contribution of this work extends beyond methodological innovation to serve as a decision-support tool for advancing real-world energy systems and policies. One limitation of SOMIT is that, although the use of medians reduces the influence of extreme outliers, the method still depends on the quality and reliability of the input data (be it subjective assessment or ACM data). Future work could enhance SOMIT by incorporating group decision-making mechanisms and by extending the method to better accommodate uncertain, imprecise, and fuzzy inputs. Beyond the two case studies presented in this work for multi-criteria renewable energy evaluation, SOMIT can be applied to a wider range of real-world energy contexts in future, such as technology portfolio selection in smart grids, location/site selection, energy storage planning, decarbonization of industry and transport, and carbon neutrality strategies. These broader applications would further demonstrate the practical benefits of SOMIT to policymakers, managers, and industry stakeholders.